\documentclass[a4paper, 12pt]{article}
\usepackage{amsfonts}
\usepackage{fancyhdr}
\usepackage{graphicx}
\usepackage{epsfig}
\usepackage{amssymb}
\usepackage{mathrsfs,amsfonts,amsmath}
\usepackage{color}

\makeatletter

\@addtoreset{equation}{section} \makeatother
\newtheorem{Theorem}{Theorem}[section]
\newtheorem{Remark}[Theorem]{Remark}
\newtheorem{Lemma}[Theorem]{Lemma}
\newtheorem{Assumption}[Theorem]{Assumption}
\newtheorem{Definition}[Theorem]{Definition}

\begin{document}
\title{\textbf{Polynomial stability of exact solution and a numerical method for stochastic differential equations with time-dependent delay
\footnote{Supported by Natural Science Foundation of China (NSFC
11601025).}}}
\author{ Guangqiang Lan\footnote{Corresponding author: Email:
langq@mail.buct.edu.cn.},
 Fang Xia,\ Qiushi Wang
\\ \small School of Science, Beijing University of Chemical Technology, Beijing 100029, China}

\date{}

\maketitle

\begin{abstract}
Polynomial stability of exact solution and modified truncated Euler-Maruyama method for stochastic differential equations with time-dependent delay are investigated in this paper. By using the well known discrete semimartingale convergence theorem, sufficient conditions are obtained for both bounded and unbounded delay $\delta$ to ensure the polynomial stability of the corresponding numerical approximation. Examples are presented to illustrate the conclusion.
\end{abstract}



\section{Introduction and main result}

\noindent

Asymptotic stability of stochastic differential delay equations has attracted more and more attention in recent years \cite{KPE,LC,Mao}. Since the exact solution is usually difficult to obtain, properties of the corresponding numerical simulations become more and more hot topics. There are plenty of papers devoted to the exponential stability of the different types of numerical solutions. For example, \cite{WMS} considered the almost sure exponential stability of Euler and backward Euler methods for stochastic delay differential equations, \cite{CZ} investigated exponential mean-square stability of two-step Maruyama methods for stochastic delay differential equations, \cite{QH} studied delay-dependent exponential stability of the backward Euler method for nonlinear stochastic delay differential equations. One can also refer to other literatures for exponential stability of numerical solutions, see e.g. \cite{CW,LY,MY,ZW} and references therein. When there is no delay, \cite{LFM} obtained the polynomial stability of the classical and backward Euler method under given conditions. However, as far as we know, there are few papers concerning about the polynomial stability of the numerical solution for the underlying stochastic differential equations with unbounded delay except \cite{M}.

Recently, Mao \cite{Mao1} introduced truncated EM method for stochastic differential equation without delay, and then he obtained sufficient conditions for the strong convergence rate of it in \cite{Mao2}. Motivated by these two works, we have introduced in \cite{LX} a new numerical simulation (which we called modified truncated Euler-Maruyama method) and obtained the strong convergence rate of it. Then we investigated $p$-th moment exponential stability of it in \cite{LX1}.

In this paper, we will first extend modified truncated Euler-Maruyama method for stochastic differential equations to that of stochastic differential equations with time dependent delay (both bounded and unbounded cases), and then we will investigate the almost sure and mean square polynomial stability of the given modified truncated Euler-Maruyama method.

Let $(\Omega,\mathscr{F},P)$ be a complete probability space with a filtration $\left\{\mathscr F_t\right\}_{t \ge 0}$ satisfying the usual conditions (i.e. it is right continuous and $\mathscr{F}_0$ contains all $\mathbb{P}$-null sets). Let $\tau\ge0$ be a constant and $C\left(\left[-\tau,0\right];\mathbb{R}^n\right)$ the space of all continuous functions from $[-\tau ,0]$ to $\mathbb{R}{^n}$ with the norm $\left\| \phi  \right\| = {\sup _{ - \tau  \le \theta  \le 0}}\left| {\phi \left( \theta  \right)} \right|$. Denote by $C_{{{\mathscr F}_0}}^b\left( {[ - \tau ,0];\mathbb{R}{^n}} \right)$ the family of bounded, $\mathscr {F}_0$ measurable, $C\left( {\left[ { - \tau ,0} \right];\mathbb{R}{^n}} \right)$-valued random variables. Let $B\left( t \right)$ be a $d$-dimensional standard Brownian motion.

Consider the following stochastic differential delay equations:
\begin{equation}\label{sdde}
dx(t)=f(x(t),x(t-\delta(t)),t)dt+g(x(t),x(t-\delta(t)),t)dB(t), t \ge 0,
\end{equation}
with the initial value
\[x_0=\xi = \left\{\xi(\theta),\theta\in[-\tau,0]\right\}\in C_{\mathscr F_0}^b\left([-\tau,0];\mathbb{R}^n\right),\]
where $\delta(t)\in C^1(\mathbb{R}_+,\mathbb{R}_+)$ such that $\delta(0)=\tau$, $\mathbb{E}||\xi||^2<\infty,$ moreover, $f:\mathbb{R}{^n} \times \mathbb{R}{^n} \times [0, + \infty ) \to \mathbb{R}{^n}$ and $g:\mathbb{R}{^n} \times \mathbb{R}{^n} \times [0, + \infty ) \to \mathbb{R}{^n}\otimes\mathbb{R}{^d}$ are Borel measurable vector and matrix valued functions, respectively.

Notice that stochastic pantograph equation is a special case of the above stochastic delay differential equation (\ref{sdde}) with unbounded memory (i.e. $\delta(t)=t-qt$, $0<q<1$ and $x_0=x(0)\in\mathbb{R}^n$ is a $\mathscr F_0$ measurable random variable).

We always assume that
\[f(0,0,t)\equiv0,\quad g(0,0,t)\equiv0,\]
which implies that $X\equiv0$ is the trivial solution of equation (\ref{sdde}). And we assume that
\begin{equation}\label{delta}\delta'(t)\le\eta<1.\end{equation}
This implies that $t-\delta(t)$ is strictly increasing on $[0,\infty)$.

We impose two standing hypotheses on $f$ and $g$ in this paper.

\begin{Assumption}\label{A1}
The coefficients $f$ and $g$ satisfy local Lipschitz condition for any fixed $t>0$, that is, for each
$R$ and $t$ there is $L_{R,t}>0$ such that
\begin{equation}\label{local}|f(x,y,t)-f(x',y',t)|
\vee|g(x,y,t)-g(x',y',t)|\le L_{R,t}\left(|x-x'|+|y-y'|\right)\end{equation} for all
$|x|\vee|x'|\vee|y|\vee|y'|\le R.$

\end{Assumption}
Here the norm of a matrix $A$ is denoted by $|A|=\sqrt{\textrm{trace}(A^\textrm{T}A)}$.

It is obvious that for any fixed $t$, ${L_{R,t}}$ is an increasing function with respect to $R$.

\begin{Assumption}\label{A2}
There exist positive constants $K$ and $\lambda_i, i=1,2,3$ such that
\begin{equation}\label{tiaojian1}
2\langle x,f(x,y,t)\rangle+|g(x,y,t)|^2\le\frac{K(1+t)^{-\lambda_0}-\lambda_1|x|^2+\lambda_2|y|^2}{1+t}
\end{equation}
for any $t \ge 0$ and $x ,y \in\mathbb{R}{^n}$.
\end{Assumption}

Now let us give the modified truncated Euler-Maruyama method for equation (\ref{sdde}).

Without loss of generality, for the given $\tau\ge0$, we can choose the step size $\Delta \in {(0,1)}$ suitably such that there exists a positive integer $m$ such that $\tau=m\Delta$.

For $\Delta^*>0$, let $h(\Delta)$ be a strictly positive decreasing function $h:(0,\Delta^*]\to(0,\infty)$ such that
\begin{equation}\label{tiaojian}
\lim_{\Delta\to0}h(\Delta)=\infty.
\end{equation}

We now define $f_\Delta$ for any $\Delta > 0$
\begin{equation}\label{dingyi}
{f_\Delta }(x,y,t) = \left\{ {\begin{array}{*{20}{l}}
{f(x,y,t),\left| x \right| \vee \left| y \right| \le h\left( \Delta  \right)},\\
{\frac{{\left| x \right| \vee \left| y \right|}}{{h\left( \Delta  \right)}}f(\frac{{h\left( \Delta  \right)}}{{\left| x \right| \vee \left| y \right|}}x,\frac{{h\left( \Delta  \right)}}{{\left| x \right| \vee \left| y \right|}}y,t),\left| x \right| \vee \left| y \right| > h\left( \Delta  \right).}
\end{array}} \right.
\end{equation}
$g_\Delta$ is defined in the same way as $f_\Delta$ .

Now, we can define the modified truncated EM (MTEM) method $X_k^\Delta\approx x(k\Delta)$ by setting  $X_k^\Delta=\xi(k\Delta)$ for every integer $k=-m,\cdots,0$ ,
and
\begin{equation}\label{num1}
X_{k+1}^\Delta=X_k^\Delta+ {f_\Delta }\left( {X_k^\Delta ,X_{k-[\frac{\delta(k\Delta)}{\Delta}]}^\Delta ,k\Delta } \right)\Delta  + {g_\Delta }\left( {X_k^\Delta ,X_{k -[\frac{\delta(k\Delta)}{\Delta}]}^\Delta ,k\Delta } \right)\Delta {B_k}
\end{equation}
for every integer $k=1,2,\cdots$, where $[x]$ is the integer part of $x$ and $\Delta {B_k}=B((k+1)\Delta)-B(k\Delta)$.

\begin{Definition}
The solution $x(t,\xi)$ to equation (\ref{sdde}) is said to be $p$-th moment exponentially stable if there exists $\gamma>0$ such that
$$\limsup_{t\rightarrow\infty}\frac{\log \mathbb{E}|x(t)|^p}{\log (1+t)}<-\gamma$$
for any initial value $x_0=\xi\in
C_{\mathscr{F}_0}^b([-\tau,0];\mathbb{R}^n).$

If $p=2,$ it is said mean square polynomially stable. It is said to be almost surely polynomially stable if for almost all $\omega\in\Omega,$
$$\limsup_{t\rightarrow\infty}\frac{\log |x(t)|}{\log (1+t)}<-\gamma.$$
\end{Definition}

\begin{Definition}
We say that the MTEM approximation $X_k^\Delta$ (\ref{num1}) is $p$-th moment polynomially stable if there exist $\Delta^*>0$ and $\gamma>0$ such that for any $0<\Delta\le\Delta^*$
\begin{equation}\label{jixian}\limsup_{k\rightarrow\infty}\frac{\log\mathbb{E}(|X_k^\Delta|^p)}{\log(1+k\Delta)}<-\gamma.\end{equation}

If $p=2,$ we say that $X_k^\Delta$ (\ref{num1}) is mean square polynomially stable. It is said to be almost surely polynomially stable if (\ref{jixian}) is replaced by
\begin{equation}\label{jixian1}\limsup_{k\rightarrow\infty}\frac{\log|X_k^\Delta|}{\log(1+k\Delta)}<-\gamma, \quad a.s..\end{equation}
\end{Definition}

Notice that under Assumption \ref{A1}, it follows that for any initial value, there exists a unique maximal local solution to equation (\ref{sdde}). Indeed, we have the following

\begin{Theorem}\label{jingque}
Assume that (\ref{delta}), \textbf{Assumption1.1} and \textbf{Assumption1.2} hold with $\lambda _1-\frac{\lambda_2}{1-\eta}>0$. Then for any initial condition $\xi\in C_{\mathscr{F}_0}^b([-\tau,0];\mathbb{R}^n)$ there exists a unique global solution $x=\{x(t),t\in[-\tau,\infty)\}$ of equation (\ref{sdde}).

Moreover, if $\delta\le\tau$, or $\delta$ is unbounded, then the solution $x(t)$ is also almost surely and mean square polynomially stable. That is, for any $\gamma\in(0,\gamma^*),$
\begin{equation}\label{jixian11}\limsup_{t\rightarrow\infty}\frac{\log |x(t)|}{\log(1+t)}\le-\frac{\gamma}{2}\ a.s., \ \textrm{and}\ \limsup_{t\rightarrow\infty}\frac{\log\mathbb{E}(|x(t)|^2)}{\log(1+t)}\le-\gamma, \end{equation}
where $\gamma^*=\lambda_0\wedge\gamma_0$ if $\delta\le\tau$, here $\gamma_0$ is the unique positive solution to $\gamma_0-\lambda_1
+\frac{\lambda_2\left(1\vee(1+\tau)^{\gamma_0-1}\right)}{1-\eta}=0,$ and $\gamma^*=\lambda_0\wedge\left(\lambda _1-\frac{\lambda_2}{1-\eta}\right)\wedge1$ if $\delta$ is unbounded.
\end{Theorem}

So a natural question raises: \textit{Does the MTEM method $X_k^\Delta$ (\ref{num1}) replicates the polynomial stability of the equation (\ref{sdde}) under given conditions?} The answer is YES.

Now we are ready to present our first main result about MTEM method (\ref{num1}). Suppose $\tau>0$ and $\delta(t)\le\tau, \forall t\ge0.$

\begin{Theorem}\label{stability}
Assume that (\ref{delta}), \textbf{Assumption1.1} and \textbf{Assumption1.2} hold with $\lambda _1-\lambda_2([(1-\eta)^{-1}]+1)>0$, and the local Lipshitz constant $L_{R,t}$ satisfies
\begin{equation}\label{l}\lim_{R\to\infty}\sup_{t\ge0}(1+t)L^2_{R,t}h^{-1}(R)=0.\end{equation}

Then the MTEM approximation (\ref{num1}) is both almost surely and mean square polynomially stable. Precisely, for any $0<\varepsilon<\frac{\lambda_1-\lambda_2([(1-\eta)^{-1}]+1)}{[(1-\eta)^{-1}]+2}$, there exists $\Delta^*>0$ and $\tilde{C}>0$ such that for any $\Delta\in (0,\Delta^*]$ and $C\in(0,\tilde{C})$,
\begin{equation}\label{zuizhong1}
\limsup_{k\to\infty}\frac{\log
|X_{k}^\Delta|}{\log(1+k\Delta)}\le-\frac{C}{2}, a.s.,\ \textrm{and}\ \limsup_{k\to\infty}\frac{\log\mathbb{E}
(|X_{k}^\Delta|^2)}{\log(1+k\Delta)}\le-C,
\end{equation}
where $\tilde{C}=\tilde{C}_0\wedge\lambda_0$, and $\tilde{C}_0$ is the unique positive solution of the following equation
$$\tilde{C}_0-(\lambda _1-\varepsilon)+(\lambda_2+\varepsilon)([(1-\eta)^{-1}]+1)(1+\tau)^{\tilde{C}_0}=0.$$
\end{Theorem}

The second main result about MTEM method (\ref{num1}) is for unbounded $\delta$.

\begin{Theorem}\label{stability1}
Let all assumptions in Theorem \ref{stability} hold. If $\delta(t)$ is unbounded, then the MTEM approximation (\ref{num1}) is both almost surely and mean square polynomially stable. That is, for any $$\max\left\{0,\frac{\lambda_1-\lambda_2([(1-\eta)^{-1}]+1)-1}{[(1-\eta)^{-1}]+2}\right\}
<\varepsilon<\frac{\lambda_1-\lambda_2([(1-\eta)^{-1}]+1)}{[(1-\eta)^{-1}]+2},$$ there exists $\Delta^*>0$ and $\tilde{C}>0$ such that for any $\Delta\in (0,\Delta^*]$ and $C\in(0,\tilde{C})$,
\begin{equation}\label{zuizhong2}
\limsup_{k\to\infty}\frac{\log
|X_{k}^\Delta|}{\log(1+k\Delta)}\le-\frac{C}{2}, a.s.,\ \textrm{and}\ \limsup_{k\to\infty}\frac{\log\mathbb{E}
(|X_{k}^\Delta|^2)}{\log(1+k\Delta)}\le-C,
\end{equation}
where $\tilde{C}=\tilde{C}_0\wedge\lambda_0$, $\tilde{C}_0=\lambda_1-\varepsilon-(\lambda_2+\varepsilon)([(1-\eta)^{-1}]+1)(<1).$
\end{Theorem}

\begin{Remark}
When there is no delay term, (\ref{tiaojian1}) becomes to $2\langle x,f(x,t)\rangle+|g(x,t)|^2\le K(1+t)^{-(\lambda_0+1)}-\lambda_1(1+t)^{-1}|x|^2.$ Then (2.5), (2.6) of condition 2.3 in \cite{LFM} implies this special cases if we take $K_1=\lambda_1=(\lambda_0+1)/2,$ and we do not need the linear growth condition (2.4) there. So our results cover that of \cite{LFM}. We also remark that Theorem \ref{stability1} can not cover Theorem \ref{stability} since in Theorem \ref{stability} the rate of polynomial stability could be larger than $1$ while in Theorem \ref{stability1} it must be smaller than $1$. Moreover, since $\lambda _1-\lambda_2([(1-\eta)^{-1}]+1)>0$ implies $\lambda _1-\frac{\lambda_2}{1-\eta}>0,$ then under the same assumptions of Theorem \ref{stability} or \ref{stability1}, there exists a unique global solution $x(t)$ to equation (\ref{sdde}) and the solution $x(t)$ is also almost surely and mean square polynomially stable. So the MTEM method $X_k^\Delta$ (\ref{num1}) replicates the polynomial stability of the equation (\ref{sdde}) under given conditions.
\end{Remark}

The rest of the paper is organized as follows. In Section 2, Theorem \ref{jingque} will be proved. Section 3 gives some lemmas which will play important roles in the proof of Theorem \ref{stability} and \ref{stability1}. We will then prove in Section 4 the almost sure and mean-square polynomial stability of the given numerical approximation when $\delta\le\tau$. Section 5 deals with the unbounded $\delta$. The last section gives numerical examples and simulations to illustrate the conclusion.

\section{The existence, uniqueness and polynomial stability of the exact solution}
Before proving Theorem \ref{jingque}, we give the following continuous semimartingale convergence theorem established in \cite{LS}, which is critical in the proof of almost sure polynomial stability of the exact solution.

\begin{Lemma}\label{l0}
Let $A(t), U(t)$ be two continuous $\mathscr{F}_{t}$ adapted increasing processes on $t\ge0$ with $A(0)=U(0)=0$ a.s. Let $M(t)$ be a real-valued continuous local martingale with $M(0)=0$ a.s. Let $\xi$ be a nonnegative $\mathscr{F}_{0}$-measurable random variable. Assume that $\{X(t)\}$ is a nonnegative semimartingale with the Doob-Meyer decomposition
$$X(t)=\xi+A(t)-U(t)+M(t), t\ge0.$$

If $\lim_{i\to\infty}A(t)<\infty$ a.s., then
$$\lim_{t\to\infty}X(t)<\infty\ \textrm{and}\ \lim_{t\to\infty}U(t),\quad a.s.$$
\end{Lemma}

\textbf{Proof of Theorem \ref{jingque}} Since local Lipschitz condition holds, it follows that for any initial condition $\xi\in C_{\mathscr{F}_0}^b([-\delta(0),0];\mathbb{R}^n)$ there exists a unique maximal local solution $x(t)$ on $[-\delta(0),\tau_e)$, where $\tau_e$ is the explosion time. Let $k_0>0$ be sufficiently large such that $||\xi||=\sup_{\theta\in[-\delta(0),0]}|\xi(\theta)|\le k_0.$ For each integer $k\ge k_0,$ define the stopping time
$$\tau_k=\inf\{t\in[0,\tau_e):|x(t)|\ge k\},\quad \inf \emptyset=\infty.$$

It is obvious that $\tau_k$ is increasing as $k\to\infty$. Then $\tau_\infty:=\lim_{k\to\infty}\tau_k$ exists and $\tau_\infty\le\tau_e$ a.s. So we only need to prove $\tau_\infty=\infty$ which implies that the solution $x(t),[-\delta(0),\infty)$ does not explode in finite time.

By It\^o formula and Assumption \ref{A2}, for any $k\ge k_0,$ it follows that
$$\aligned |x(t\wedge\tau_k)|^2&=|x(0)|^2+\int_0^{t\wedge\tau_k}\left(2\langle x(s),f(x(s),x(s-\delta(s)),s)+|g(x(s),x(s-\delta(s)),s)|^2\right)ds\\&
\quad+2\int_0^{t\wedge\tau_k}\langle x(s),g(x(s),x(s-\delta(s)),s)dB(s)\rangle\\&
\le |x(0)|^2+\int_0^{t\wedge\tau_k}\left(\frac{K}{(1+s)^{1+\lambda_0}}-\frac{\lambda_1}{1+s}|x(s)|^2+\frac{\lambda_2}{1+s}|x(s-\delta(s))|^2\right)ds\\&
\quad+2\int_0^{t\wedge\tau_k}\langle x(s),g(x(s),x(s-\delta(s)),s)dB(s)\rangle.\endaligned$$

Then taking expectation on both sides, we obtain for every $t\ge0,$
$$\aligned\mathbb{E}|x(t\wedge\tau_k)|^2\le \mathbb{E}||\xi||^2+\frac{K}{\lambda_0}+\mathbb{E}\int_0^{t\wedge\tau_k}\left(-\frac{\lambda_1}{1+s}|x(s)|^2+\frac{\lambda_2}{1+s}|x(s-\delta(s))|^2\right)ds\endaligned$$

Since $\delta'(t)\le\eta<1,$ then there exists a unique positive $t_0$ such that $t_0=\delta(t_0)$. Thus,
\begin{equation}\label{jie0}\aligned\mathbb{E}\int_0^{t\wedge\tau_k}\frac{\lambda_2}{1+s}|x(s-\delta(s))|^2ds&
=\mathbb{E}\int_{t_0}^{t\wedge\tau_k}\frac{\lambda_2}{1+s}|x(s-\delta(s))|^2ds\\&
\quad+\mathbb{E}\int_0^{t_0}\frac{\lambda_2}{1+s}|x(s-\delta(s))|^2ds\\&
\le \mathbb{E}\int_0^{t\wedge\tau_k}\frac{\lambda_2}{1+r+\delta(s)}|x(r)|^2\frac{1}{1-\delta'(s)}dr+\lambda_2\mathbb{E}||\xi||^2t_0\\&
\le\frac{\lambda_2}{1-\eta}\mathbb{E}\int_0^{t\wedge\tau_k}\frac{1}{1+r}|x(r)|^2dr+\lambda_2\mathbb{E}||\xi||^2t_0.\endaligned\end{equation}

Therefore, $\lambda_1-\frac{\lambda_2}{1-\eta}>0$ yields
\begin{equation}\label{jie}\aligned\mathbb{E}|x(t\wedge\tau_k)|^2&\le \mathbb{E}||\xi||^2+\frac{K}{\lambda_0}+\lambda_2\mathbb{E}||\xi||^2t_0-\left(\lambda_1-\frac{\lambda_2}{1-\eta}\right)\mathbb{E}\int_0^{t\wedge\tau_k}\frac{1}{1+s}|x(s)|^2ds\\&
\le \mathbb{E}||\xi||^2(1+\lambda_2t_0)+\frac{K}{\lambda_0}.\endaligned\end{equation}

Letting $k\to\infty$ we have
\begin{equation}\label{jie1}\aligned\mathbb{E}|x(t\wedge\tau_\infty)|^2&
\le \mathbb{E}||\xi||^2(1+\lambda_2t_0)+\frac{K}{\lambda_0}.\endaligned\end{equation}

Now if $P(\tau_\infty<\infty)>0$ then for some $T>0, P(\tau_\infty<T)>0.$

Taking $t=T$ in (\ref{jie1}), we have
\begin{equation}\aligned P(\tau_\infty<T)|x(\tau_\infty)|^2\le\mathbb{E}|x(\tau_\infty)1_{\tau_\infty<T}|^2\le \mathbb{E}|x(t\wedge\tau_\infty)|^2
\le \mathbb{E}||\xi||^2(1+\lambda_2T)<\infty,\endaligned\end{equation}
which is impossible since $x(\tau_\infty)=\infty.$ So $P(\tau_\infty=\infty)=1,$ as required.

Now let us prove the almost sure and mean square polynomial stability.

By using It\^o formula and Assumption \ref{A2} again, for any $0<\gamma<\lambda_0$, we have
\begin{equation}\label{jq}\aligned (1+t)^\gamma|x(t)|^2&
\le||\xi||^2+\gamma\int_0^{t}(1+s)^{\gamma-1}|x(s)|^2ds+K\int_0^{t}(1+s)^{\gamma-\lambda_0-1}ds\\&
\quad+\int_0^{t}(1+s)^{\gamma}\left(-\frac{\lambda_1}{1+s}|x(s)|^2+\frac{\lambda_2}{1+s}|x(s-\delta(s))|^2\right)ds+M(t)\\&
\le||\xi||^2+\frac{K((1+t)^{\gamma-\lambda_0}-1)}{\gamma-\lambda_0}+(\gamma-\lambda_1)\int_0^{t}(1+s)^{\gamma-1}|x(s)|^2\\&
\quad+\lambda_2\int_0^t(1+s)^{\gamma-1}|x(s-\delta(s))|^2ds+M(t),\endaligned\end{equation}
where $M(t)=2\int_0^{t}(1+s)^{\gamma}\langle x(s),g(x(s),x(s-\delta(s)),s)dB(s)\rangle$ is a continuous local martingale with $M(0)=0.$

\textit{\textbf{Case 1:}} If $\delta\le\tau,$ then for any $\gamma>0,$ we have
$$\aligned\int_0^t(1+s)^{\gamma-1}|x(s-\delta(s))|^2ds&\le ||\xi||^2t_0+\int_{t_0}^t(1+s)^{\gamma-1}|x(s-\delta(s))|^2ds\\&
\le||\xi||^2t_0+\frac{1\vee(1+\tau)^{\gamma-1}}{1-\eta}\int_{0}^t(1+r)^{\gamma-1}|x(r)|^2dr.\endaligned$$
Here $a\vee b=\max\{a,b\}$.

Therefore, for any $\gamma\in(0,\lambda_0)$
$$\aligned (1+t)^\gamma|x(t)|^2&
\le||\xi||^2(1+\lambda_2t_0)+\frac{K}{\lambda_0-\gamma}\\&\quad+\left(\gamma-\lambda_1
+\frac{\lambda_2\left(1\vee(1+\tau)^{\gamma-1}\right)}{1-\eta}\right)\int_0^{t}(1+s)^{\gamma-1}|x(s)|^2ds+M(t).\endaligned$$

Since $\lambda _1-\frac{\lambda_2}{1-\eta}>0,$ then there exists a unique $\gamma_0>0$ such that for any $\gamma\in(0,\gamma_0)$,
$$\gamma-\lambda_1
+\frac{\lambda_2\left(1\vee(1+\tau)^{\gamma-1}\right)}{1-\eta}\le0.$$

Consequently, for any $\gamma\in(0,\lambda_0\wedge\gamma_0)$
$$(1+t)^{\gamma}|x(t)|^2\le ||\xi||^2(1+\lambda_2t_0)+\frac{K}{\lambda_0-\gamma}+M(t).$$

\textit{\textbf{Case 2:}} If $\delta$ is unbounded, then for any $\gamma<1$, similar to (\ref{jie0}), we have
$$\aligned\lambda_2\int_0^t(1+s)^{\gamma-1}|x(s-\delta(s))|^2ds&\le \lambda_2||\xi||^2t_0+\lambda_2\int_{t_0}^t(1+s)^{\gamma-1}|x(s-\delta(s))|^2ds\\&
\le \lambda_2||\xi||^2t_0+\frac{\lambda_2}{1-\eta}\int_{0}^t(1+r)^{\gamma-1}|x(r)|^2dr.\endaligned$$

Thus for any $\gamma<1\wedge\lambda_0$,
$$\aligned (1+t)^\gamma|x(t)|^2&
\le||\xi||^2(1+\lambda_2t_0)+\frac{K}{\lambda_0-\gamma}+\left(\gamma-\lambda_1+\frac{\lambda_2}{1-\eta}\right)\int_0^{t}(1+s)^{\gamma-1}|x(s)|^2ds.\endaligned$$

Since $\lambda_1-\frac{\lambda_2}{1-\eta}>0,$ then for any $\gamma<(-\lambda_1+\frac{\lambda_2}{1-\eta})\wedge1\wedge\lambda_0,$ we have
$$\aligned (1+t)^{\gamma}|x(t)|^2&
\le||\xi||^2(1+\lambda_2t_0)+\frac{K}{\lambda_0-\gamma}+M(t).\endaligned$$

So by Lemma \ref{l0}, the exact solution $x(t)$ is almost surely polynomially stable in both cases. For mean square polynomial stability, we only need to take expectation on both sides of the above equation, then
$$\aligned (1+t)^{\gamma}\mathbb{E}|x(t)|^2&
\le\mathbb{E}||\xi||^2(1+\lambda_2t_0)+\frac{K}{\lambda_0-\gamma}.\endaligned$$
We complete the proof. $\square$

\section{Some useful lemmas}

\noindent
To prove our main results, let us present some useful lemmas.

We first introduce the so called discrete semimartingale convergence theorem (cf. \cite{Mao,WMS}), which is essential in proving the main results in this paper.
\begin{Lemma}\label{l4}
Let $\{A_i\}, \{U_i\}$ be two sequences of nonnegative random variables such that both $A_i$ and $U_i$ are $\mathscr{F}_{i-1}$-measurable for $i=1,2,\cdots,$ and $A_0=U_0=0$ a.s. Let $M_i$ be a real-valued local martingale with $M_0=0$ a.s. Let $\xi$ be a nonnegative $\mathscr{F}_{0}$-measurable random variable. Assume that $\{X_i\}$ is a nonnegative semimartingale with the Doob-Meyer decomposition
$$X_i=\xi+A_i-U_i+M_i.$$

If $\lim_{i\to\infty}A_i<\infty$ a.s., then for almost all $\omega\in\Omega,$
$$\lim_{i\to\infty}X_i<\infty\qquad\lim_{i\to\infty}U_i<\infty,$$
that is, both $X_i$ and $U_i$ converge to finite random variables.
\end{Lemma}

Now let introduce the following lemma, which reveals completely the significance of the constant $\eta$ introduced by (\ref{delta})

For any given $\tau>0$, we can choose $\Delta>0$ such that $\frac{\tau}{\Delta}=m$, where $m$ is a positive integer. Then $$i-\left[\frac{\delta(i\Delta)}{\Delta}\right]\ge i-\frac{\delta(i\Delta)}{\Delta}\ge-\frac{\delta(0)}{\Delta}\ge-\frac{\tau}{\Delta}= -m.$$

We have

\begin{Lemma}\label{l3}
Suppose (\ref{delta}) holds. For an arbitrary but fixed $i\in\{0,1,2,\cdots\}$, let $i-\left[\frac{\delta(i\Delta)}{\Delta}\right]=a,$ where $a\in\{-m,-m+1,\cdots,0,1,\cdots,i\}$ Then
\begin{equation}\label{geshu}
\#\left\{j\in\{0,1,2,\cdots\}:j-\left[\frac{\delta(j\Delta)}{\Delta}\right]=a\right\}\le[(1-\eta)^{-1}]+1,
\end{equation}
where $\#S$ denotes the number of elements of the set $S$.
\end{Lemma}

The proof of Lemma \ref{l3} can be found in \cite{M1}.

For the modified truncated function $f_\Delta$ and $g_\Delta$, we have the following global Lipschitz continuity.

\begin{Lemma}\label{L1}
Suppose the local Lipshitz condition (\ref{local}) holds. Then for any fixed $\Delta>0$,
\begin{equation}\label{global}
|f_\Delta(x,y,t)- f_\Delta(\bar{x},\bar{y},t)|\le 5L_{h(\Delta ),t}(|x- \bar{x}|+|y-\bar{y}|).
\end{equation}
\end{Lemma}

\textbf{Proof}\ For any $x, y, \bar{x}, \bar{y}\in\mathbb{R}^d$, there are three cases: $(|x|\vee|y|)\vee(|\bar{x}|\vee|\bar{y}|)\le h(\Delta)$, $(|x|\vee|y|)\wedge(|\bar{x}|\vee|\bar{y}|)> h(\Delta)$ and one of $|x|\vee|y|$ and $|\bar{x}|\vee|\bar{y}|$ is no greater than $h(\Delta)$ and the other of them is greater than $h(\Delta)$.

If $|x|\vee|y|\vee|\bar{x}|\vee|\bar{y}|\le h(\Delta)$, then (\ref{global}) holds naturally by (\ref{local}).

Now assume $(|x|\vee|y|)\wedge(|\bar{x}|\vee|\bar{y}|)> h(\Delta)$.

Since

$$\left|\frac{h(\Delta)}{|x|\vee|y|}x\right|\vee\left|\frac{h(\Delta)}{|x|\vee|y|}y\right|
\vee\left|\frac{h(\Delta)}{|\bar{x}|\vee|\bar{y}|}\bar{x}\right|\vee\left|\frac{h(\Delta)}{|\bar{x}|\vee|\bar{y}|}\bar{y}\right|\le h(\Delta),$$
then by (\ref{local}), we have

$$\aligned |f_\Delta(x,y,t)-f_\Delta(\bar{x},\bar{y},t)|&=\left|\frac{|x|\vee|y|}{h(\Delta)} f\left(\frac{h(\Delta)}{|x|\vee|y|}(x,y),t\right)-\frac{|\bar{x}|\vee|\bar{y}|}{h(\Delta)} f\left(\frac{h(\Delta)}{|\bar{x}|\vee|\bar{y}|}(\bar{x},\bar{y}),t\right)\right|\\&
\le \frac{|x|\vee|y|}{h(\Delta)}\left|f\left(\frac{h(\Delta)}{|x|\vee|y|}(x,y),t\right)
-f\left(\frac{h(\Delta)}{|\bar{x}|\vee|\bar{y}|}(\bar{x},\bar{y}),t\right)\right|\\&
\quad+ \left|f\left(\frac{h(\Delta)}{|\bar{x}|\vee|\bar{y}|}(\bar{x},\bar{y}),t\right)\right|\cdot
\left|\frac{|x|\vee|y|-|\bar{x}|\vee|\bar{y}|}{h(\Delta)}\right|\\&
\le \frac{|x|\vee|y|}{h(\Delta)}\cdot L_{h(\Delta),t}\left(\left|\frac{h(\Delta)}{|x|\vee|y|}x-\frac{h(\Delta)}{|\bar{x}|\vee|\bar{y}|}\bar{x}\right|\right.\\&
\qquad\qquad\qquad\qquad+\left.\left|\frac{h(\Delta)}{|x|\vee|y|}y-\frac{h(\Delta)}{|\bar{x}|\vee|\bar{y}|}\bar{y}\right|\right)\\&
\quad+ 2L_{h(\Delta),t}h(\Delta)\left|\frac{|x|\vee|y|-|\bar{x}|\vee|\bar{y}|}{h(\Delta)}\right|\\&
= L_{h(\Delta),t}\left(\left|x-\frac{|x|\vee|y|}{|\bar{x}|\vee|\bar{y}|}\bar{x}\right|
+\left|y-\frac{|x|\vee|y|}{|\bar{x}|\vee|\bar{y}|}\bar{y}\right|\right)\\&
\quad+ 2L_{h(\Delta),t}\left||x|\vee|y|-|\bar{x}|\vee|\bar{y}|\right|\\&
\le L_{h(\Delta),t}\left(\left|x-\bar{x}\right|+|y-\bar{y}|+2\left||x|\vee|y|-|\bar{x}|\vee|\bar{y}|\right|\right)
\\&
\quad+ 2L_{h(\Delta),t}\left||x|\vee|y|-|\bar{x}|\vee|\bar{y}|\right|.\endaligned$$
Here and from now on, $(a(x,y),t):=(ax,ay,t)$.

Since
$$\aligned\left||x|\vee|y|-|\bar{x}|\vee|\bar{y}|\right|
&=\left|\frac{|x|+|y|+||x|-|y||}{2}-\frac{|\bar{x}|+|\bar{y}|+||\bar{x}|-|\bar{y}||}{2}\right|\\&
\le\frac{1}{2}(|x-\bar{x}|+|y-\bar{y}|+||x|-|y|-|\bar{x}|+|\bar{y}||)\\&
\le |x-\bar{x}|+|y-\bar{y}|,\endaligned$$
then
$$|f_\Delta(x,y,t)-f_\Delta(\bar{x},\bar{y},t)|\le5L_{h(\Delta),t}(|x-\bar{x}+|y-\bar{y}|).$$

Finally, without loss of generality, suppose that $|x|\vee|y|\le h(\Delta)<|\bar{x}|\vee|\bar{y}|.$ Then we have
$$\aligned |f_\Delta(x,y,t)-f_\Delta(\bar{x},\bar{y},t)|&=\left|f(x,y,t)-\frac{|\bar{x}|\vee|\bar{y}|}{h(\Delta)} f\left(\frac{h(\Delta)}{|\bar{x}|\vee|\bar{y}|}(\bar{x},\bar{y}),t\right)\right|\\&
\le \left|f(x,y,t)- f\left(\frac{h(\Delta)}{|\bar{x}|\vee|\bar{y}|}(\bar{x},\bar{y}),t\right)\right|\\&
\quad+\left|f\left(\frac{h(\Delta)}{|\bar{x}|\vee|\bar{y}|}(\bar{x},\bar{y}),t\right)\right|
\left|1-\frac{|\bar{x}|\vee|\bar{y}|}{h(\Delta)}\right|\\&
\le L_{h(\Delta),t}\left(\left|x-\frac{h(\Delta)}{|\bar{x}|\vee|\bar{y}|}\bar{x}\right|
+\left|y-\frac{h(\Delta)}{|\bar{x}|\vee|\bar{y}|}\bar{y}\right|\right)\\&
\quad+ 2h(\Delta)L_{h(\Delta),t}\left|1-\frac{|\bar{x}|\vee|\bar{y}|}{h(\Delta)}\right|\\&
\le L_{h(\Delta),t}\left(|x-\bar{x}|+|\bar{x}|\left|1-\frac{h(\Delta)}{|\bar{x}|\vee|\bar{y}|}\right|\right.\\&
\quad\left.+|y-\bar{y}|+|y|\left|1-\frac{h(\Delta)}{|\bar{x}|\vee|\bar{y}|}\right|\right)\\&
\quad+2L_{h(\Delta),t}\left|h(\Delta)-|\bar{x}|\vee|\bar{y}|\right|\\&
\le L_{h(\Delta),t}\left(|x-\bar{x}|+|y-\bar{y}|+4|h(\Delta)-|\bar{x}|\vee|\bar{y}||\right).
\endaligned$$

Since $|x|\vee|y|\le h(\Delta)<|\bar{x}|\vee|\bar{y}|$, then $$|h(\Delta)-|\bar{x}|\vee|\bar{y}||\le\left||x|\vee|y|-|\bar{x}|\vee|\bar{y}|\right|\le|x-\bar{x}|+|y-\bar{y}|$$

Therefore,
$$\aligned|f_\Delta(x,y,t)-f_\Delta(\bar{x},\bar{y},t)|&\le 5L_{h(\Delta),t}(|x-\bar{x}+|y-\bar{y}|).\endaligned$$

Similarly, we can prove that $g_\Delta$ is globally Lipschitz continuous with the same Lipschitz constant $5L_{h(\Delta),t}$.
 We complete the proof. $\square$

 \begin{Lemma}\label{L2}
Suppose (\ref{tiaojian1}) holds. Then for any fixed $\Delta>0$,
\begin{equation}\label{tiaojian2}2\langle x,f_\Delta(x,y,t)\rangle+|g_\Delta(x,y,t)|^2\le\frac{K(1+t)^{-\lambda_0}-(\lambda_1-\frac{K}{h^2(\Delta)})|x|^2+(\lambda_2+\frac{K}{h^2(\Delta)})|y|^2}{1+t} \end{equation}
\end{Lemma}

\textbf{Proof}\ On the one hand, (\ref{tiaojian2}) holds naturally by (\ref{tiaojian1}) and the definitions of $f_\Delta$ and $g_\Delta$ if $|x|\vee|y|\le h(\Delta)$.

On the other hand, if $|x|\vee|y|>h(\Delta)$, then
\[\begin{array}{*{20}{l}}
2\langle x,f_\Delta(x,y,t)\rangle+|g_\Delta(x,y,t)|^2 &=2\left\langle x,\frac{|x|\vee|y|}{h(\Delta)}f\left(\frac{h(\Delta)}{|x|\vee|y|}x,\frac{h(\Delta)}{|x|\vee|y|}y,t\right)\right\rangle\\
&\quad+\frac{|x|^2\vee|y|^2}{h^2(\Delta)}\left|g\left(\frac{h(\Delta)}{|x|\vee|y|}x,\frac{h(\Delta)}{|x|\vee|y|}y,t\right)\right|^2\\
&=2\langle x,\frac{1}{a}f(ax,ay,t)\rangle+\frac{1}{a^2}|g(ax,ay,t)|^2\\
&=\frac{1}{a^2}(2\langle ax,f(ax,ay,t)\rangle+|g(ax,ay,t)|^2)
\end{array}\]
where $a=\frac{h(\Delta)}{|x|\vee|y|}$.
Then by using (\ref{tiaojian1}), it follows that
\[\begin{array}{*{20}{l}}
2\langle x,f_\Delta(x,y,t)\rangle+|g_\Delta(x,y,t)|^2&\le\frac{1}{a^2}{(1+t)^{-1}}\left(K(1+t)^{-\lambda_0}-\lambda_1|ax|^2+\lambda_2|ay|^2 \right)\\
&\le\frac{K(|x|^2+|y|^2)}{h^2(\Delta)(1+t)}+(1+t)^{-1}(-\lambda_1|x|^2+\lambda_2|y|^2)\\&
=\frac{-\left(\lambda_1-\frac{K}{h^2(\Delta)}\right)|x|^2+\left(\lambda_2+\frac{K}{h^2(\Delta)}\right)|y|^2}{1+t},
\end{array}\]
as required. $\square$

\section{Polynomial stability of $X_k^\Delta$ when $\delta$ is bounded}

\textbf{Proof of Theorem \ref{stability}:} By the definition of MTEM (\ref{num1}), we have
\begin{equation}\label{diedai}
\aligned|X_{k + 1}^\Delta|^2&=|X_k^\Delta|^2+|f_{\Delta,k}|^2\Delta ^2 + 2\langle X_k^\Delta ,f_{\Delta,k}\Delta \rangle +|g_{\Delta,k}|^2\Delta + M_k,\endaligned
\end{equation}
where
\[f_{\Delta,k}={{f_\Delta }\left( {X_k^\Delta ,X_{k -\left[\frac{\delta(k\Delta)}{\Delta}\right] }^\Delta ,k\Delta } \right)}, \]
\[g_{\Delta,k}={{g_\Delta }\left( {X_k^\Delta ,X_{k -\left[\frac{\delta(k\Delta)}{\Delta}\right] }^\Delta ,k\Delta } \right)}, \]
and
\[\begin{array}{*{20}{l}}
{m_k}:= &{2\left\langle {X_k^\Delta  + f_{\Delta,k }\Delta,g_{\Delta,k }\Delta {B_k}} \right\rangle}+(|g_{\Delta,k }\Delta {B_k} |^2- |g_{\Delta,k }|^2\Delta).
\end{array} \]

Then by using (\ref{global}) and (\ref{tiaojian2}), we have
\[\aligned
|X_{k + 1}^\Delta|^2& =|X_k^\Delta|^2 +(2\langle X_k^\Delta,f_{\Delta,k}\rangle+|g_{\Delta,k}|^2)\Delta+|f_{\Delta,k}|^2\Delta^2+m_k\\&
\le |X_k^\Delta|^2+\frac{K}{(1+k\Delta)^{\lambda_0+1}}+\left(\frac{-(\lambda_1-\frac{K}{h^2(\Delta)})\Delta}{1+k\Delta}+50L_{h(\Delta),k\Delta}^2\Delta^2\right)| X_k^\Delta|^2 \\&\quad+ \left(\frac{(\lambda_2+\frac{K}{h^2(\Delta)})\Delta}{1+k\Delta}+50L_{h(\Delta),k\Delta}^2\Delta ^2 \right)\left|X^\Delta_{k-\left[\frac{\delta(k\Delta)}{\Delta}\right]}\right|^2+m_k.
\endaligned\]

By (\ref{l}), if we set $R=h(\Delta)$ and $t=k\Delta$, then for $\Delta \to 0$ (thus $R\to\infty$),
\[(1+k\Delta)L_{h(\Delta),k\Delta}^2\Delta=(1+t)L^2_{R,t}h^{-1}(R)\to 0.\]

That is $L_{h(\Delta),k\Delta}^2\Delta=o(\frac{1}{1+k\Delta})$ for any fixed $k$. Since $h(\Delta)\to\infty$ as $\Delta\to0,$ then for any sufficiently small $\varepsilon>0$, there exists a $\Delta^*\in(0,1)$ small enough such that for all $\Delta\in (0,\Delta^*)$,$\left(50(1+k\Delta)L_{h(\Delta),k\Delta}^2\Delta\right)\vee\frac{K}{h^2(\Delta)}\le\varepsilon$.

Thus,
\[|X_{k+1}^\Delta|^2\le|X_k^\Delta|^2+\frac{K}{(1+k\Delta)^{\lambda_0+1}}-\frac{(\lambda _1-\varepsilon)\Delta }{1 + k\Delta }|X_k^\Delta|^2+ \frac{(\lambda _2+\varepsilon)\Delta }{1 + k\Delta}\left|X_{k-\left[\frac{\delta(k\Delta)}{\Delta}\right]}^\Delta\right|^2+m_k.\]

Observe that, for an arbitrary constant $C>0$, if we multiply both sides by ${\left( {1 + \left( {k + 1} \right)\Delta } \right)^C}$, then
\[\begin{array}{*{20}{l}}
(1+(k+1)\Delta)^C|X_{k+1}^\Delta|^2&{\le {\left( {1 + k\Delta } \right)^C}{\left| {X_k^\Delta } \right|^2}+\frac{K(1+(k+1)\Delta)^C}{(1+k\Delta)^{\lambda_0+1}}+ C_k{\left| {X_k^\Delta } \right|^2}}\\
&{\quad+ {\left( {1 + \left( {k + 1} \right)\Delta } \right)^C}\frac{{\left( {{\lambda _2} + \varepsilon } \right)\Delta }}{{1 + k\Delta }}{\left| {X_{k - \left[\frac{\delta(k\Delta)}{\Delta}\right]}^\Delta } \right|^2+{\left( {1 + \left( {k + 1} \right)\Delta } \right)^C}m_k,}}
\end{array}\]
where $C_k= {{{\left( {1 + \left( {k + 1} \right)\Delta } \right)}^C} - {{\left( {1 + k\Delta } \right)}^C} - {{\left( {1 + \left( {k + 1} \right)\Delta } \right)}^C}\frac{{\left( {{\lambda _1}-\varepsilon }\right)\Delta }}{{1 + k\Delta }}}.$

Thus,
\begin{equation}\label{ditui}\aligned
(1 + k\Delta)^C|X_k^\Delta|^2&\le|X_0^\Delta|^2+K\sum\limits_{i=0}^{\infty}\frac{(1+(i+1)\Delta)^C}{(1+i\Delta)^{\lambda_0+1}}+\sum\limits_{i = 0}^{k - 1} C_i|X_i^\Delta|^2 \\&
\quad+\sum\limits_{i = 0}^{k-1}(1+(i+1)\Delta)^C\frac{(\lambda _2+\varepsilon)\Delta}{1+i\Delta}\left| X^\Delta_{i-\left[\frac{\delta(i\Delta)}{\Delta}\right]}\right|^2+M_k
\endaligned\end{equation}
where $M_k=\sum_{i=0}^{k-1}(1 + (i+1)\Delta)^Cm_i.$

It is obvious that for any $\Delta>0,$ $\{M_k,\mathscr{F}_{k\Delta}\}_{k\ge0}$ is a local martingale with $M_0=0$. Moreover, if $\lambda_0>C$ then $\sum\limits_{i=0}^{\infty}\frac{(1+(i+1)\Delta)^C}{(1+i\Delta)^{\lambda_0+1}}<\infty$

Notice that by (\ref{delta}), there exists unique $t_0\ge0$ such that $t_0=\delta(t_0)$. Then for any fixed $\Delta$ there exists a unique $i_0$ (independent of $k$) such that $\forall 0\le i\le i_0, j\ge i_0+1,$
$$i-\left[\frac{\delta(i\Delta)}{\Delta}\right]< 0\ \textrm{and}\ j-\left[\frac{\delta(j\Delta)}{\Delta}\right]\ge 0.$$

Thus
\begin{equation}\label{d}\aligned &\quad
 \sum\limits_{i = 0}^{k-1} {{{\left( {1 + (i+1)\Delta } \right)}^C}} \frac{{\left( {{\lambda _2} + \varepsilon } \right)\Delta }}{{1 + i\Delta }}\mathbb{E}{{\left| {X_{i-\left[\frac{\delta(i\Delta)}{\Delta}\right]}^\Delta } \right|}^2}\\&
 = \sum\limits_{i\ : i-\left[\frac{\delta(i\Delta)}{\Delta}\right]< 0}{{{\left( {1 + (i+1)\Delta } \right)}^C}} \frac{{\left( {{\lambda _2} + \varepsilon } \right)\Delta }}{{1 + i\Delta }}\mathbb{E}{{\left| {X_{i-\left[\frac{\delta(i\Delta)}{\Delta}\right]}^\Delta } \right|}^2}\\
&\quad + \sum\limits_{{i\ : i-\left[\frac{\delta(i\Delta)}{\Delta}\right]\ge 0}} {{{\left( {1 + (i+1)\Delta } \right)}^C}} \frac{{\left( {{\lambda _2} + \varepsilon } \right)\Delta }}{{1 + i\Delta }}\mathbb{E}{{\left| {X_{i-\left[\frac{\delta(i\Delta)}{\Delta}\right]}^\Delta } \right|}^2}\\
& =: D_1 + D_2.\endaligned\end{equation}

Notice that
$$\aligned D_1&
\le \sum\limits_{i\le i_0}{{{\left( {1 + (i+1)\Delta } \right)}^C}} \frac{{\left( {{\lambda _2} + \varepsilon } \right)\Delta }}{{1 + i\Delta }}{{\left| {X_{i-\left[\frac{\delta(i\Delta)}{\Delta}\right]}^\Delta } \right|}^2}\\&
\le (1+(i_0+1)\Delta)^C(\lambda_2+\varepsilon)\Delta\sum\limits_{i\le i_0}{{\left| {X_{i-\left[\frac{\delta(i\Delta)}{\Delta}\right]}^\Delta } \right|}^2}\\&
=:D'_1.\endaligned$$
It is obvious that $D'_1$ is a nonnegative $\mathscr{F}_0$ measurable random variable. And
$$D_2=\sum\limits_{{i\ : i-\left[\frac{\delta(i\Delta)}{\Delta}\right]\ge 0}} {{{\left( {1 + \left( x_i +y_i+1 \right)\Delta } \right)}^C}} \frac{{\left({{\lambda _2} + \varepsilon } \right)\Delta }}{{\left(1 +x_i+y_i \right)\Delta }}{{\left| {X_{x_i}^\Delta } \right|}^2,}$$
where $x_i:=i- \left[\frac{\delta(i\Delta)}{\Delta}\right]$ and $y_i:=\left[\frac{\delta(i\Delta)}{\Delta}\right]$ are both nonnegative.

Moreover, since for any $x,y\ge0$
$$\frac{\frac{(1+(x+1)\Delta+y)^C}{1+x\Delta+y}}{\frac{(1+(x+1)\Delta)^C}{1+x\Delta}}
=\left(\frac{1+(x+1)\Delta+y}{1+(x+1)\Delta}\right)^C\times\frac{1+x\Delta}{1+x\Delta+y}\le(1+y)^C,$$
then
$$\aligned&\quad\frac{\left(1+\left(x_i+y_i+1\right)\Delta
\right)^C}{1+\left(x_i+y_i\right)\Delta}\\&
\le \left(1+y_i\Delta\right)^C\frac{\left(1+\left(x_i+1\right)\Delta
\right)^C}{1+x_i\Delta}\\&
\le\left(1+\tau\right)^C\frac{\left(1+\left(x_i+1\right)\Delta
\right)^C}{1+x_i\Delta}.\endaligned$$

Therefore, Lemma \ref{l3} yields that
\begin{equation}\label{kz1}\aligned&\quad\sum\limits_{i = 0}^{k-1}(1+(i+1)\Delta)^C\frac{(\lambda_2+\varepsilon)\Delta}{1+i\Delta}\left| X^\Delta_{i-\left[\frac{\delta(i\Delta)}{\Delta}\right]}\right|^2\\&
\le D'_1+(\lambda _2+\varepsilon)\Delta(1+\tau)^C\sum\limits_{i = i_0+1}^{k-1}
\frac{\left(1+\left(x_i+1\right)\Delta
\right)^C}{1+x_i\Delta
}\left| X^\Delta_{x_i}\right|^2\\&
\le D'_1+(\lambda _2+\varepsilon)\Delta(1+\tau)^C([(1-\eta)^{-1}]+1)\sum\limits_{j=0}^{k-1}
\frac{\left(1+\left(j+1\right)\Delta
\right)^C}{1+j\Delta
}\left| X^\Delta_{j}\right|^2.\endaligned\end{equation}

Now by (\ref{ditui}), it follows that
\begin{equation}\label{ditui1}\aligned
(1 + k\Delta)^C|X_k^\Delta|^2&\le|X_0^\Delta|^2+K\sum\limits_{i=0}^{\infty}\frac{(1+(i+1)\Delta)^C}{(1+i\Delta)^{\lambda_0+1}}+D'_1+\sum\limits_{i = 0}^{k - 1} K_i|X_i^\Delta|^2+M_k\\&
\le Y+\sum\limits_{i = 0}^{k - 1} K_i|X_i^\Delta|^2+M_k,
\endaligned\end{equation}
where $Y:=|X_0^\Delta|^2+K\sum\limits_{i=0}^{\infty}\frac{(1+(i+1)\Delta)^C}{(1+i\Delta)^{\lambda_0+1}}+D'_1$ is also a nonnegative $\mathscr{F}_0$ measurable random variable and
$$\aligned K_i&=C_i+(\lambda _2+\varepsilon)\Delta(1+\tau)^C([(1-\eta)^{-1}]+1)\frac{\left(1+\left(i+1\right)\Delta
\right)^C}{1+i\Delta}\\&
=(1+(i+1)\Delta)^C\left[1-\left(1+\frac{\Delta}{1+i\Delta}\right)^{-C}-\frac{(\lambda _1-\varepsilon)\Delta }{1 + i\Delta}\right]\\&
\quad+(\lambda _2+\varepsilon)\Delta(1+\tau)^C([(1-\eta)^{-1}]+1)\frac{\left(1+\left(i+1\right)\Delta
\right)^C}{1+i\Delta}\\&
\le (1+(i+1)\Delta)^C \left[1-\left(1-\frac{C\Delta}{1+i\Delta}\right)-\frac{(\lambda _1-\varepsilon)\Delta }{1 + i\Delta}\right]\\&
\quad+(\lambda _2+\varepsilon)\Delta(1+\tau)^C([(1-\eta)^{-1}]+1)\frac{\left(1+\left(i+1\right)\Delta
\right)^C}{1+i\Delta}\\&
=\frac{(1+(i+1)\Delta)^C}{1+i\Delta}\left[C-(\lambda _1-\varepsilon)+(\lambda_2+\varepsilon)(1+\tau)^C([(1-\eta)^{-1}]+1)\right]\Delta.\endaligned$$

Notice that since $\lambda _1-\lambda_2([(1-\eta)^{-1}]+1)>0,$ then we can choose $\varepsilon>0$ sufficiently small such that
$$\lambda _1-\varepsilon-(\lambda_2+\varepsilon)([(1-\eta)^{-1}]+1)>0.$$

On the other hand, since for any fixed $a>b>0,$ $f(x):=x-a+b(1+\tau)^x$ is strictly increasing on $[0,\infty]$, $f(0)=0$ and $\lim_{x\to\infty}f(x)=\infty$, then there exists a unique $x_0>0$ such that $f(x_0)=b-a<0$. Let $a=\lambda _1-\varepsilon$ and $b=(\lambda_2+\varepsilon)([(1-\eta)^{-1}]+1)$. Then there exists a unique $\tilde{C}_0>0$ such that
$$\tilde{C}_0-(\lambda _1-\varepsilon)+(\lambda_2+\varepsilon)([(1-\eta)^{-1}]+1)(1+\tau)^{\tilde{C}_0}=0.$$

Thus for any $C\in(0,\lambda_0\wedge\tilde{C}_0), $(\ref{ditui}) implies
\begin{equation}\label{ditui2}\aligned
(1 + k\Delta)^C|X_k^\Delta|^2&\le Y+M_k.
\endaligned\end{equation}

Consequently, the well known discrete semimartingale convergence theorem (Lemma \ref{l4}) yields that, for almost all $\omega\in\Omega,$
\[\limsup_{k\to\infty}(1+k\Delta)^C|X_k^\Delta|^2<\infty,\]
which implies (\ref{zuizhong1}).

For mean square polynomial stability, taking expectation on both sides of (\ref{ditui2}) yields that
\begin{equation}\label{ditui3}\aligned
(1+k\Delta)^C\mathbb{E}|X_k^\Delta|^2&\le \mathbb{E}Y.
\endaligned\end{equation}

Since
$$\aligned \mathbb{E}Y&=\mathbb{E}\left(|X_0^\Delta|^2+K\sum\limits_{i=0}^{\infty}\frac{(1+(i+1)\Delta)^C}{(1+i\Delta)^{\lambda_0+1}}+D'_1\right)\\&
\le||\xi||^2+K\sum\limits_{i=0}^{\infty}\frac{(1+(i+1)\Delta)^C}{(1+i\Delta)^{\lambda_0+1}}+(i_0+1)(1+(i_0+1)\Delta)^C(\lambda_2+\varepsilon)\Delta||\xi||^2\\&
<\infty,\endaligned$$
and it is obvious that $\mathbb{E}Y$ is independent of $k$, then
\begin{equation}\aligned
\limsup_{k\to\infty}(1+k\Delta)^C\mathbb{E}|X_k^\Delta|^2\le\mathbb{E}Y<\infty,
\endaligned\end{equation}
as required. $\square$

\section{Polynomial stability of $X_k^\Delta$ when $\delta$ is unbounded}

\textbf{Proof of Theorem \ref{stability1}:} The idea is same as the proof of Theorem \ref{stability}. We will first obtain an inequality similar to (\ref{ditui2}), then use the discrete semimartingale convergence theorem.

Even if $\delta$ is unbounded now, by repeating the proof of Theorem \ref{stability} word by word, we have for any $C>0$
\begin{equation}\label{ditui4}\aligned
(1 + k\Delta)^C|X_k^\Delta|^2&\le|X_0^\Delta|^2+K\sum\limits_{i=0}^{\infty}\frac{(1+(i+1)\Delta)^C}{(1+i\Delta)^{\lambda_0+1}}+\sum\limits_{i = 0}^{k - 1} C_i|X_i^\Delta|^2 \\&
\quad+\sum\limits_{i = 0}^{k-1}(1+(i+1)\Delta)^C\frac{(\lambda _2+\varepsilon)\Delta}{1+i\Delta}\left| X^\Delta_{i-\left[\frac{\delta(i\Delta)}{\Delta}\right]}\right|^2+M_k
\endaligned\end{equation}
where $C_i= {{{\left( {1 + \left( {i + 1} \right)\Delta } \right)}^C} - {{\left( {1 + i\Delta } \right)}^C} - {{\left( {1 + \left( {i + 1} \right)\Delta } \right)}^C}\frac{{\left( {{\lambda _1}-\varepsilon }\right)\Delta }}{{1 + i\Delta }}}$, and $M_k=\sum_{i=0}^{k-1}(1 + (i+1)\Delta)^Cm_i$ is a $\{\mathscr{F}_{k\Delta}\}_{k\ge0}$ martingale and $M_0=0.$

Notice that
$$\aligned C_i&= (1+(i+1)\Delta)^C\left(1-\left(1+\frac{\Delta}{1+i\Delta}\right)^{-C}-\frac{(\lambda _1-\varepsilon)\Delta }{1 + i\Delta}\right)\\&
\le (1+(i+1)\Delta)^C \left(1-\left(1-\frac{C\Delta}{1+i\Delta}\right)-\frac{(\lambda _1-\varepsilon)\Delta }{1 + i\Delta}\right)\\&
=\frac{(1+(i+1)\Delta)^C}{1+i\Delta}\left(C-(\lambda_1-\varepsilon)\right)\Delta.\endaligned$$

Then
\begin{equation}\label{kz}\aligned
(1 + k\Delta)^C|X_k^\Delta|^2&\le |X_0^\Delta|^2+K\sum\limits_{i=0}^{k-1}\frac{(1+(i+1)\Delta)^C}{(1+i\Delta)^{\lambda_0+1}}\\&
\quad+\left(C-(\lambda_1-\varepsilon)\right)\Delta\sum\limits_{i = 0}^{k - 1} \frac{(1+(i+1)\Delta)^C}{1+i\Delta}|X_i^\Delta|^2 \\&
\quad+\sum\limits_{i = 0}^{k-1}(1+(i+1)\Delta)^C\frac{(\lambda _2+\varepsilon)\Delta}{1+i\Delta}\left| X^\Delta_{i-\left[\frac{\delta(i\Delta)}{\Delta}\right]}\right|^2+M_k.
\endaligned\end{equation}

Since $\delta$ is unbounded in this case, then the inequality (\ref{kz1}) does not hold any more. Thus we have to estimate
$$\sum\limits_{i = 0}^{k-1}(1+(i+1)\Delta)^C\frac{(\lambda _2+\varepsilon)\Delta}{1+i\Delta}\left| X^\Delta_{i-\left[\frac{\delta(i\Delta)}{\Delta}\right]}\right|^2$$
in another way.

Notice that (\ref{d}) still holds in this case, and
$$\aligned D_1&
\le \sum\limits_{i\le i_0}{{{\left( {1 + (i+1)\Delta } \right)}^C}} \frac{{\left( {{\lambda _2} + \varepsilon } \right)\Delta }}{{1 + i\Delta }}{{\left| {X_{i-\left[\frac{\delta(i\Delta)}{\Delta}\right]}^\Delta } \right|}^2}\\&
\le (1+(i_0+1)\Delta)^C(\lambda_2+\varepsilon)\Delta\sum\limits_{i\le i_0}{{\left| {X_{i-\left[\frac{\delta(i\Delta)}{\Delta}\right]}^\Delta } \right|}^2}\\&
=:D'_1.\endaligned$$
It is obvious that $D'_1$ is a $\mathscr{F}_0$ measurable random variable.

On the other hand,
$$D_2=\sum\limits_{{i\ : i-\left[\frac{\delta(i\Delta)}{\Delta}\right]\ge 0}} {{{\left( {1 + \left( x_i +y_i+1 \right)\Delta } \right)}^C}} \frac{{\left({{\lambda _2} + \varepsilon } \right)\Delta }}{{\left(1 +x_i+y_i \right)\Delta }}{{\left| {X_{x_i}^\Delta } \right|}^2,}$$
where $x_i:=i- \left[\frac{\delta(i\Delta)}{\Delta}\right]$ and $y_i:=\left[\frac{\delta(i\Delta)}{\Delta}\right]$.

Bearing in mind that
$$\left(\frac{(1+(x+1)\Delta)^C}{1+x\Delta}\right)'=\frac{(1+(1+x)\Delta)^{C-1}(C-1+x\Delta(C-1)-\Delta)\Delta}{(1+x\Delta)^2}.$$

If $C\le1,$ then
$$\left(\frac{(1+(x+1)\Delta)^C}{1+x\Delta}\right)'\le0,\ \forall x\ge0.$$

Thus for $C\le1$, it follows that
\[\frac{\left(1 + \left( x_i +y_i+1 \right)\Delta\right)^C}{1 +\left(x_i+y_i \right)\Delta}\le\frac{(1+(x_i+1)\Delta)^C}{1+x_i\Delta}.\]

By Lemma \ref{l3}, it follows that
$$\aligned D_2&\le\left(\lambda_2+\varepsilon\right)\Delta\sum\limits_{{i\ : i-\left[\frac{\delta(i\Delta)}{\Delta}\right]\ge 0}} \frac{\left( 1 + \left( x_i+1 \right)\Delta\right)^C}{\left(1 +x_i\right)\Delta }\left|X_{x_i}^\Delta\right|^2\\&
\le ([(1-\eta)^{-1}]+1)\left(\lambda_2+\varepsilon\right)\Delta\sum_{i=0}^{k-1}\frac{\left(1+\left(i+1\right)\Delta\right)^C}{\left(1 +i\right)\Delta}\left|X_{i}^\Delta\right|^2,\endaligned$$

Then (\ref{kz}) implies that
\begin{equation}\label{kz2}(1+k\Delta)^C|X_k^\Delta|^2\le|X_0^\Delta|^2
+K\sum\limits_{i=0}^{\infty}\frac{(1+(i+1)\Delta)^C}{(1+i\Delta)^{\lambda_0+1}}+D'_1+\sum_{i = 0}^{k-1}K_i{\left|X_i^\Delta\right|^2}+M_k,\end{equation}
where
\begin{equation}\label{KI}
\aligned K_i&:=\frac{\left(1+\left(i+1\right)\Delta\right)^C}{\left(1 +i\right)\Delta}\left[C-(\lambda_1-\varepsilon)+ ([(1-\eta)^{-1}]+1)\left(\lambda_2+\varepsilon\right)\right]\Delta.\endaligned
\end{equation}

Now if $0<\lambda_1-\lambda_2([(1-\eta)^{-1}]+1)\le1,$ then for any $\varepsilon>0$ small enough ($<\frac{\lambda_1-\lambda_2([(1-\eta)^{-1}]+1)}{[(1-\eta)^{-1}]+2}$), we have
$$0<(\lambda_1-\varepsilon)- ([(1-\eta)^{-1}]+1)\left(\lambda_2+\varepsilon\right)<\lambda_1-\lambda_2([(1-\eta)^{-1}]+1)\le1.$$

If $\lambda_1-\lambda_2([(1-\eta)^{-1}]+1)>1,$ then we can choose $\frac{\lambda_1-\lambda_2([(1-\eta)^{-1}]+1)-1}{[(1-\eta)^{-1}]+2}<\varepsilon
<\frac{\lambda_1-\lambda_2([(1-\eta)^{-1}]+1)}{[(1-\eta)^{-1}]+2}$.

Thus we still have
$$0<(\lambda_1-\varepsilon)- ([(1-\eta)^{-1}]+1)\left(\lambda_2+\varepsilon\right)<1.$$

So if $0<\lambda_1-\lambda_2([(1-\eta)^{-1}]+1),$ then for any $\varepsilon>0$ such that $$\max\left\{0,\frac{\lambda_1-\lambda_2([(1-\eta)^{-1}]+1)-1}{[(1-\eta)^{-1}]+2}\right\}
<\varepsilon<\frac{\lambda_1-\lambda_2([(1-\eta)^{-1}]+1)}{[(1-\eta)^{-1}]+2},$$
it follows that
$$0<(\lambda_1-\varepsilon)- ([(1-\eta)^{-1}]+1)\left(\lambda_2+\varepsilon\right)<1.$$

Then for any
$$C\in\left(0,\left((\lambda_1-\varepsilon)- ([(1-\eta)^{-1}]+1)\left(\lambda_2+\varepsilon\right)\right)\wedge\lambda_0\right),$$
$K_i\le0$ for any $i=0,1,\cdots, k-1.$

Then we have
\[(1+k\Delta)^C|X_k^\Delta|^2\le|X_0^\Delta|^2+K\sum\limits_{i=0}^{\infty}\frac{(1+(i+1)\Delta)^C}{(1+i\Delta)^{\lambda_0+1}}+D'_1+M_k,\]

The following is same as that of Section 3. We complete the proof. $\square$

\section{Examples}

Now let us present some examples to interpret our conclusion.

\textbf{Example 1} Let $n=1, \tau>0.$ Consider the following scalar SDDE:
 \begin{equation}\label{sde1}\aligned dx(t)&=\frac{-2x(t)+\frac{1}{2}x(t-\delta(t))-x^3(t)-x(t)x^4(t-\delta(t))}{1+t}dt\\&
 \quad+\sqrt{\frac{2x^2(t)x^4(t-\delta(t))+\frac{1}{2}x^2(t-\delta(t))+2x^4(t)}{1+t}}dB_t\endaligned\end{equation}
with initial value $x_0=\{\xi(\theta), \theta\in[-\tau,0]\}\in C([-\tau,0],\mathbb{R}^n)$ and $\delta(t)=\tau+\frac{1}{2}-\frac{1}{2}e^{-t}$
In this case $f(x,y,t)=\frac{-2x+\frac{1}{2}y-x^3-xy^4}{1+t}$, $g=\sqrt{\frac{2x^2y^4+\frac{1}{2}y^2+2x^4}{1+t}}$.

It is obvious that $\delta(t)\le\tau$ and $\delta'(t)=\frac{1}{2}e^{-t}\in(0,\frac{1}{2}]$ for any $t\ge0$. So condition (\ref{delta}) holds for this $\delta$ and $\eta=\frac{1}{2}$.

Moreover,
$$\aligned2\langle x,f(x,y,t) \rangle+|g(x,y,t)|^2&=\frac{-4x^2+xy+\frac{1}{2}y^2}{1+t}\\&
\le(1 + t)^{-1}\left(-\frac{7}{2}x^2+y^2\right).\endaligned$$

That is, Assumption \ref{A2} holds for $\lambda_1=\frac{7}{2}$ and $\lambda_2=1$ and any $\lambda_0>0 (K=0).$

On the other hand, by mean value theorem for two dimensional function, we have
$$|f(x,y,t)-f(x',y',t)|\le\frac{5+4R^4}{1+t}(|x-x'|+|y-y'|)$$
and
$$|g(x,y,t)-g(x',y',t)|\le\frac{5+5R^3}{1+t}(|x-x'|+|y-y'|)$$
for all $R>0$ and $|x|\vee|x'|\vee|y|\vee|y'|\le R.$
Thus Assumption \ref{A1} holds for $L_{R,t}=\frac{5(R^4+2)}{1+t}.$

Choose $h(\Delta)=\Delta^{-\frac{1}{9}}$. Then $0<h(\Delta)\to\infty$ as $\Delta\to0,$ and $h^{-1}(R)=R^{-9}$. Thus,
$$\sup_{t\ge0}(1+t)L^2_{R,t}h^{-1}(R)=\sup_{t\ge0}\frac{25(R^4+2)^2}{(1+t)R^9}\to0$$
as $R\to\infty.$ That is, (\ref{l}) holds for such defined $h$.

Notice that in this case
$$\lambda _1-\lambda_2([(1-\eta)^{-1}]+1)=\frac{3}{2}-1\times([(1-\frac{1}{2})^{-1}]+1)=\frac{1}{2}>0.$$

Then by Theorem \ref{stability}, for any fixed $0<\varepsilon<\frac{1}{8}$, there exists $\Delta^*$ small enough such that for any $\Delta<\Delta^*$ (of course $\frac{1}{\Delta}$ is an integer) the MTEM method $X_k^\Delta$ for (\ref{sde1}) is almost surely polynomially stable with rate $\frac{\tilde{C}_0}{2}$, and it is also mean square polynomially stable with rate $\tilde{C}_0$, where $\tilde{C}_0$ is the unique positive solution to
$$\tilde{C}_0-\left(\frac{3}{2}-\varepsilon\right)+3(1+\varepsilon)(1+\tau)^{\tilde{C}_0}=0.$$

On the other hand, by Theorem \ref{jingque}, it follows that the exact solution to (\ref{sdde}) is also almost surely and mean square polynomially stable. Thus, the MTEM method replicates the polynomial stability of the the exact solution for the given SDDE (\ref{sde1}).

If we choose $\tau=1, x_0(\theta)\equiv2,$ $\theta\in[-1,0],$ $\Delta=0.1$ and $h(\Delta)=\Delta^{-\frac{1}{9}}$, then computer simulation (Matlab) for the first $5000$ steps of discrete MTEM (\ref{num1}) confirms the almost sure polynomial stability and $\frac{\log|X_k^\Delta|}{\log(1+k\Delta)}$ is less than $-1$ for $k$ large enough ($k=5000$). Notice that
$$1-\left(\frac{3}{2}-\varepsilon\right)+6(1+\varepsilon)>0$$
and
$$-\left(\frac{3}{2}-\varepsilon\right)+3(1+\varepsilon)<0.$$

Then we know that $\tilde{C}_0\in(0,1)$. Therefore $\frac{\log|X_k^\Delta|}{\log(1+k\Delta)}\le\tilde{C}_0$.

However, the numerical approximation $X_k^\Delta$ is not exponentially stable since $\frac{\log|X_k^\Delta|}{k\Delta}\to0.$

\textbf{Example 2} Let $\frac{1}{2}\le q<1$. Consider the following stochastic pantograph equation:
 \begin{equation}\label{sde2}\aligned dx(t)&=\frac{-2x(t)+\frac{1}{2}x(qt)-x^3(t)-x(t)x^4(qt)}{1+t}dt\\&
 \quad+\sqrt{\frac{2x^2(t)x^4(qt)+\frac{1}{2}x^2(qt)+2x^4(t)}{1+t}}dB_t\endaligned\end{equation}
with initial value $x(0)\in\mathbb{R}^1.$

Here $\delta(t)=t-qt$, $f(x,y,t)=\frac{-2x+\frac{1}{2}y-x^3-xy^4}{1+t}$ and $g=\sqrt{\frac{2x^2y^4+\frac{1}{2}y^2+2x^4}{1+t}}$.

It is obvious that $\delta$ is unbounded and satisfies (\ref{delta}) with $\eta=1-q$ in this case. Moreover, by Example 1, we have known that $f$ and $g$ satisfy both Assumption \ref{A1} and Assumption \ref{A2}, and (\ref{l}) holds for $h(\Delta)=\Delta^{-\frac{1}{9}}$.

Moreover, since $\eta=1-q\in(0,\frac{1}{2}],$ then
$$\lambda _1-\lambda_2([(1-\eta)^{-1}]+1)=\lambda _1-\lambda_2([q^{-1}]+1)=\frac{1}{2}>0.$$

Then by Theorem \ref{stability1}, for any
$\varepsilon\in(0,
\frac{1}{8}),$
there exists $\Delta^*>0$ such that for any $\Delta\in (0,\Delta^*]$,
\begin{equation}
\limsup_{k\to\infty}\frac{\log
|X_{k}^\Delta|}{\log(1+k\Delta)}\le-\frac{\tilde{C}_0}{2}, a.s.,
\end{equation}
and
\begin{equation}
\limsup_{k\to\infty}\frac{\log\mathbb{E}
(|X_{k}^\Delta|^2)}{\log(1+k\Delta)}\le-\tilde{C}_0,
\end{equation}
where $\tilde{C}_0=\frac{1}{2}-4\varepsilon(<1).$

If we choose $x_0(\theta)=3,$ $q=0.5,$ $\Delta=0.05$ and $h(\Delta)=\Delta^{-\frac{1}{9}}$, then computer simulation (Matlab) for the first $1000$ steps of discrete MTEM (\ref{num1}) indicates the almost sure polynomial stability and $X_k^\Delta$ is polynomially stable with $\frac{\log|X_k^\Delta|}{\log(1+k\Delta)}$ less than $-1<\frac{1}{2}-4\varepsilon=\tilde{C}_0$.

\end{document}